\documentclass[12pt]{article}
\usepackage{amssymb,amsmath}
\newtheorem{theorem}{Theorem}
\newtheorem{apptheorem}{Theorem}[section]
\newtheorem{applemma}{Lemma}[section]
\newcommand{\bull}{\mbox{$\;\;\;$\vrule height .9ex width .8ex depth -.1ex}}
\newenvironment{proof}{\par\smallbreak\noindent{\bf Proof.~}}
{\unskip\nobreak\hfill \bull \par\medbreak}
\newcommand{\shift}[1]{\mathit{shift}(#1)}
\newcommand{\fix}[1]{\mathit{fix}(#1)}
\newcommand{\function}[2]{:#1 \rightarrow #2}
\newcommand{\reals}{\mathbb{R}}
\newcommand{\setdef}[2]{\left\{ \hspace{0.5mm} #1 : \hspace{0.5mm} #2 \right\}}

\title{How Much Work Does It Take\\ To Straighten a Plane Graph Out?}

\author{M.~Kang%
\thanks{Institut f\"ur Informatik, Humboldt Universit\"at zu Berlin, 
D-10099 Berlin}\quad%
M.~Schacht%
\thanks{Institut f\"ur Informatik, Humboldt Universit\"at zu Berlin, 
D-10099 Berlin}\quad%
O.~Verbitsky%
\thanks{IAPMM, Lviv 79060, Ukraine.
Supported by an Alexander von Humboldt return fellowship.}
}

\date{}

\begin{document}

\maketitle

\begin{abstract}
We prove that if one wants to make a plane graph drawing straight-line 
then in the worst case one has to move almost all vertices.

{\it
The second vesion of this e-print includes literally the first version.
In addition, Appendix \ref{app:bound} gives an explicit bound on the
number of fixed vertices and Appendix \ref{app:overview} gives 
an overview of related work.

The final version appears as \cite{KPRSV09}.
}
\end{abstract}

We use the standard concepts of a \emph{plane graph} and
a \emph{plane embedding} (or \emph{drawing}) of
an abstract planar graph (see, e.g., \cite{Die}).
Given a plane graph $G$, we want to redraw it making all its 
edges straight line segments while keeping as many vertices on the spot 
as possible. Let $\shift G$ denote the smallest $s$ such that
we can do the job by shifting only $s$ vertices. We define $s(n)$
to be the maximum $\shift G$ over all $G$ with $n$ vertices.

The function $s(n)$ can have another interpretation closely related
to a nice web puzzle called \emph{Planarity Game} \cite{Tan}.
At the start of the game, a player sees a straight line drawing of 
a planar graph with many edge crossings. In a move s/he is allowed 
to shift one vertex to a new position; the incident edges are redrawn 
correspondingly (being all the time straight line segments). 
The objective is to obtain a crossing-free drawing. Thus,
$s(n)$ is equal to the number of moves that the player, playing optimally,
is forced to make on an $n$-vertex game instance at the worst case.

The Wagner-F\'ary-Stein theorem (see, e.g., \cite{NRa}) says that 
every $G$ has a straight line plane embedding and immediately implies
an upper bound $s(n)\le n-3$. We here aim at proving a lower bound.

Given an abstract planar graph $G$, let $\shift{G}$ denote the maximum
$\shift{G'}$ over all plane embeddings $G'$ of $G$. Thus, we are seeking
for $G$ with large $\shift{G}$. Every 4-connected planar graph $G$
is Hamiltonian (Tutte \cite{Tut}), therefore, has a matching of size
at least $(n-1)/2$ and, therefore, $\shift G\ge(n-3)/2$. An example of planar $G$
with $3k$ vertices and $\shift G\ge 2k-8$ is shown in \cite{obf}, thereby
giving us a bound $s(n)>\frac23n-10$. We now prove a much stronger bound.

\begin{theorem}\label{thm:original}
$s(n)=n(1-o(1))$.
\end{theorem}

\begin{proof}
The vertex set of a graph $G$ will be denoted by $V(G)$.
If $X\subseteq V(G)$, then $G[X]$ denotes the subgraph induced by $G$
on~$X$.

It suffices to prove that for every $k$ and every its multiple $n$
there is an $n$-vertex $G$ with $\shift G>(1-1/k)n-k^2$. 

\medskip

{\it Construction of $G$.}

\smallskip

Let $n=k(s+k)$. Let $V(G)=\bigcup_{i=1}^{s+k}V_i$ with $|V_i|=k$
for all $i$. We will describe a plane embedding of $G$
(crossing-free, not necessary straight line).
Let each $G[V_i]$ be an arbitrary maximal planar graph.
Draw these $s+k$ fragments of $G$ so that they lie in 
the outer faces of each other (a very important condition!).
Finally, add some edges to make $G$ 3-connected.
Say, we can join each pair $G[V_i]$ and $G[V_{i+1}]$ by two non-adjacent
edges and add yet another edge between $G[V_1]$ and $G[V_{s+k}]$.

This embedding is needed only to define $G$ as an abstract graph.
Once this is done, we have to specify a ``bad'' drawing
of $G$ which is far from any straight line drawing.

\medskip

{\it ``Bad'' drawing of $G$.}

\smallskip

Let $C$ be a circle. Put each $V_i$ on $C$ at the vertices
of some regular $k$-gon. The drawing is specified.\footnote{%
There is no need to describe edges; we can suppose either that 
the drawing is straight line with edge crossings as in the Planarity Game 
or that we have an arbitrary crossing-free drawing
with edges of any shape.}

\medskip

{\it Making it straight line, crossing-free: Analysis.}

\smallskip

Let $G'$ be an arbitrary straight line, crossing-free redrawing
of $G$ in the same plane. We have to show that not many vertices
of $G'$ keep the same location as they had in $G$.
Let $V'_i$ denote the location of $V_i$ in $G'$.
Denote the complement of the outer face of $G'[V'_i]$ by $T_i$.
Since $G'[V'_i]$ is a triangulation, $T_i$ is a traingle containing 
this plane graph. 
Recall that $G$ is 3-connected. By the Whitney theorem (e.g.\ \cite{Die}), 
$G'$ is equivalent to the original (defining) plane version of $G$. That is,
either these two embeddings are obtainable from one another 
by a plane homeomorphism
or this is true after changing outer face in one of them.
By construction, the regions occupied by the $G[V_i]$'s in the original
embedding are pairwise disjoint. If we change outer face, this
is still true possibly with one exception.
It follows that all but one $T_i$'s are pairwise disjoint.
Without loss of generality suppose that the possible exception is $T_{s+k}$.

Call $V'_i$ \emph{persistent} if $i < s+k$ and $|V'_i\cap V_i|\ge2$.
Since all persistent $T_i$'s are pairwise disjoint and each of them
contains a pair of vertices of some regular $k$-gon, there can be
at most $k-1$ persistent sets. It follows that the number of moved
vertices is at least
$$
\sum_{i=1}^{s+k-1}|V'_i\setminus V_i|\ge s(k-1)=
(\frac nk-k)(k-1)=(1-\frac1k)n-k^2+k,
$$
as claimed.
\end{proof}

\subsection*{Acknowledgment}
We acknowledge an important contribution of Oleg Pikhurko, who nevertheless
declined coauthorthip of this note.

\appendix

\section{An explicit bound for the number of fixed vertices}\label{app:bound}

We now reprove Theorem \ref{thm:original}, going over our original argument 
with somewhat more care and achieving two improvements. First, we obtain
an explicit bound $s(n)\ge n-2\sqrt n-1$. Second, we show that this bound
is attained by drawings with vertices occupying any prescribed set of $n$
points in weakly convex position.

By a \emph{drawing} of a planar graph $G$ we mean an arbitrary injective map
$\pi\function{V(G)}{\reals^2}$. Given a drawing $\pi$, we suppose that each 
edge $uv$ of $G$ is drawn as the straight line segment with endpoints $\pi(u)$
and $\pi(v)$.
Due to possible edge crossings and even overlaps,
$\pi$ may not be a plane drawing of $G$. 
Hence it is natural to consider a parameter
$$
\fix{G,\pi}=\max_{\pi'}|\setdef{v\in V(G)}{\pi'(v)=\pi(v)}|
$$
where the maximum is taken over all plane straight line drawings $\pi'$ of $G$.
Note a relation to our previous notation, namely $\fix{G,\pi}=n-\shift\pi$.

We will use some elementary combinatorics of integer sequences.
A sequence identified with all its cyclic shifts
will be referred to as \emph{circular}.
Subsequences of a circular sequence $S$ will be
considered also circular sequences. Note that the set of all circular
subsequences is the same for $S$ and any its shift.
The length of a $S$ will be denoted by $|S|$.

\begin{applemma}\label{lem:circle}
Let $k,s\ge1$ and $S^{k,s}$ be the circular sequence consisting of $s$ successive
blocks of the form $12\ldots k$. 
Suppose that $S$ is a subsequence of $S^{k,s}$ with no
4-subsubsequence of the form $xyxy$, where $x\ne y$. Then $|S|<k+s$.
\end{applemma}

\begin{proof}
We proceed by the double induction on $k$ and $s$.
The base case where $k=1$ and $s$ is arbitrary is trivial.
Let $k\ge2$ and consider an $S$ with no forbidden subsequence. 
If every of the $k$ elements occurs in $S$ at most once, then
$|S|\le k$ and the claimed bound is true. Otherwise, without loss of generality 
we suppose that $S$ contains $\ell\ge2$ occurrences of $k$. 
Let $A_1,\ldots,A_\ell$ (resp.\ $B_1,\ldots,B_\ell$) 
denote the parts of $S$ (resp.\ $S^{k,s}$) between these $\ell$ elements.
Thus, $|S|=\ell+\sum_{i=1}^\ell|A_i|$. 

Denote the number of elements with at least one occurrence in $A_i$ by $k_i$.
Each element $x$ occurs in at most one of the $A_i$'s because otherwise
$S$ would contain a subsequence $xkxk$. It follows that $\sum_{i=1}^\ell k_i\le k-1$.
Note that, if we append $B_i$ with an element $k$,
it will consist of blocks $12\ldots k$. Denote the number of these blocks
by $s_i$ and notice the equality $\sum_{i=1}^\ell s_i=s$.
Since $A_i$ has no forbidden subsequence, we have $|A_i|\le k_i+s_i-1$.
If $k_i\ge1$, this follows from the induction assumption because
$A_i$ can be regarded as a subsequence of $S^{k_i,s_i}$.
If $k_i=0$, this is also true because then $|A_i|=0$.
Summarizing, we obtain $|S|\le\ell+\sum_{i=1}^\ell(k_i+s_i-1)\le\ell+(k-1)+s-\ell<k+s$.
\end{proof}

\begin{apptheorem}\label{thm:fix}
Let $k\ge3$, $n=k^2$, and $H$ be a 3-connected plane graph with $n$ vertices
having the following property: Its vertex
set can be split into $k$ equal parts $V(H)=V_1\cup\ldots\cup V_k$ so that
each $H[V_i]$ is a triangulation and these $k$ triangulations lie in 
the outer faces of each other. 
Let $X$ be an arbitrary set of $n$ points
on the boundary $\Gamma$ of a convex plane body. Then there is
a drawing $\pi\function{V(H)}X$ such that $\fix{H,\pi}\le2\sqrt n+1$.
\end{apptheorem}

\begin{proof}
Let $X=\{x_1,\ldots,x_n\}$, where the points in $X$ are numbered in the order
of their appearance along $\Gamma$. Fix $\pi$ to be an arbitrary map
such that $\pi(V_i)=\{x_i,x_{i+k},x_{i+2k},\ldots,x_{i+(k-1)k}\}$ for each $i\le k$. 

Let $\pi'$ be an arbitrary crossing-free straight line redrawing
of $H$. We have to show that not many vertices
of $H$ keep the same location in $\pi'$ as they had in $\pi$.
Denote $A_i=\setdef{\pi(v)}{v\in V_i,\ \pi(v)=\pi'(v)}$ and $A=\bigcup_{i=1}^kA_i$. 
The union $A$ consists of exactly those vertices that
keep their position under transition from $\pi$ to $\pi'$.
Thus, we have to bound the number of vertices in $A$ from above.

Denote the complement of the outer face of $H[V_i]$ in $\pi'$ by $T_i$.
Since $H[V_i]$ is a triangulation, $T_i$ is a triangle containing all $\pi'(V_i)$.
Recall that $H$ is 3-connected. By Whitney's theorem (see, e.g., \cite{Die}), 
$\pi'$ is equivalent to the original plane embedding of $H$, which we denote by $\delta$.
This means that one of the following two cases occurs:
\begin{description}
\item[\rm A]
$\pi'$ is obtainable from $\delta$ by a plane homeomorphism.
\item[\rm B]
$\pi'$ is obtainable by a plane homeomorphism from $\delta_F$,
where $F$ is an inner face of $\delta$ and $\delta_F$ is an embedding of $H$
obtained from $\delta$ by making the face $F$ outer.
\end{description}
By construction, the regions occupied by the $H[V_i]$'s are pairwise disjoint
in $\delta$. For $\pi'$ this implies that, if we have Case A, then
all $k$ triangles $T_i$ are pairwise disjoint. The same holds true
in Case B if $F$ is not a face of any $H[V_i]$-fragment.
If $F$ is a face of some $H[V_j]$-fragment, then the $T_i$'s are pairwise disjoint
with one exception for the triangle $T_j$, which contains all the others.

Consider first the case that all the triangles are pairwise disjoint.
Since $A_i\subset T_i$, the convex hulls of these sets of points are pairwise disjoint.
Label each $x_j$ by the index $i$ for which $x_j\in\pi(V_i)$
and consider the circular sequence of these labels in the order of their
appearance along $\Gamma$. This is exactly the sequence $S^{k,k}$ as in 
Lemma~\ref{lem:circle}. Let $S$ be the subsequence corresponding to
the points in $A$. Since the points in $A_i$ are labeled by $i$, we see
that $S$ has no subsequence of the form $xyxy$. By Lemma~\ref{lem:circle},
$|A|=|S|<2k$.

Consider now the case that the triangles $T_i$ are pairwise disjoint 
with the exception, say, for $T_k$.
Let $a,b,c\in V_k$ be the vertices on the boundary of the outer face 
of $H[V_k]$ in $\delta$. Let $T$ denote the geometric triangle with
vertices $\pi'(a)$, $\pi'(b)$, and $\pi'(c)$. Note that
$A_i\subset T$ for all $i<k$. Note also that $A_k\subset \Gamma\setminus T$.
The set $\Gamma\setminus T$ consists of at most three continuous components;
denote the corresponding parts of $A_k$ by $A'_k,A'_{k+1},A'_{k+2}$. 
Consider the circular sequence $S$ as above with the following
modification: the vertices in $A'_{k+1}$ are relabeled with $k+1$ and
the vertices in $A'_{k+2}$ are relabeled with $k+2$ (the vertices in $A'_k$
keep label $k$). This modification rules out any $xyxy$-subsequence.
Note that the modified $S$ can be considered a subsequence of $S^{k+2,k}$.
By Lemma~\ref{lem:circle}, we have $|A|\le2k+1$.
\end{proof}

\section{Related work}\label{app:overview}

Given a planar graph $G$, define
$
\fix G=\min_\pi\fix{G,\pi}
$,
where the minimum is taken over all drawings of $G$.
In other words, $\fix G$ is the maximum number of vertices
which can be fixed in any drawing of $G$ while ``untangling'' it.
Note that $\shift G=n-\fix G$.

The \emph{cycle} (resp.\ \emph{path}; \emph{empty graph})
on $n$ vertices will be denoted by $C_n$ (resp.\ $P_n$; $E_n$).
Recall that the \emph{join} of vertex-disjoint graphs $G$ and $H$
is the graph $G*H$ consisting of the union of $G$ and $H$ and all
edges between $V(G)$ and $V(H)$. The graphs $W_n=C_{n-1}*E_1$ (resp.\
$F_n=P_{n-1}*E_1$; $S_n=E_{n-1}*E_1$) are known as \emph{wheels} (resp.\
\emph{fans}; \emph{stars}). By $kG$ we denote
the disjoint union of $k$ copies of a graph~$G$.

Pach and Tardos~\cite{PTa} were first who established a principal fact:
Some graphs can be drawn so that, in order to untangle them, one has
to shift almost all their vertices. In fact, 
this is already true for cycles.
More precisely,  Pach and Tardos~\cite{PTa} proved that
$$
\fix{C_n}=O((n\log n)^{2/3}).
$$
This bound is nearly optimal, as shown by Cibulka~\cite{Cib}.

The best known upper bounds are of the form $\fix G=O(\sqrt n)$.
Goaoc et al.\ \cite{Goaoc} showed it for certain triangulations. 
More specifically, they proved that
$$
\fix{P_{n-2}*P_2}<\sqrt n+2.
$$

Shortly after \cite{Goaoc} and independently of it,
there appeared the first version of the current e-print.
We constructed 3-connected
planar graphs $H_n$ with $\fix{H_n}=o(n)$. Though no explicit bound was
specified in that version, a simple analysis of our construction
reveals that
$$
\fix{H_n}\le2\sqrt n+1,
$$
see Appendix \ref{app:bound}.
While $H_n$ is not as simple as $P_{n-2}*P_2$ and the subsequent examples in the
literature, the construction of $H_n$'s has the advantage that it can ensure
certain special properties of these graphs, as bounded vertex degrees.
By a later result of Cibulka \cite{Cib},
for graphs with bounded vertex degrees we have $\fix G=O(\sqrt n(\log n)^{3/2})$
whenever their diameter is logarithmic. Note in this respect that $H_n$
has diameter~$\Omega(\sqrt n)$.

In subsequent papers \cite{SWo,Bose} examples of graphs with small $\fix G$
were found in special classes of planar graphs,
as outerplanar and even acyclic graphs.
Spillner and Wolff \cite{SWo} showed for the fan graph that
\begin{equation}\label{eq:fans}
\fix{F_n}<2\sqrt n+1
\end{equation}
and Bose et al.\ \cite{Bose} established for the star forest with $n=k^2$ vertices that
\begin{equation}\label{eq:stars}
\fix{kS_k}\le3\sqrt n-3.
\end{equation}
Finally, Cibulka \cite{Cib} proved that
$$
\fix{G}=O((n\log n)^{2/3})
$$
for all 3-connected planar graphs.

Improving a result of Spillner and Wolff~\cite{SWo},
Bose et al.~\cite{Bose} showed that
$$
\fix G\ge(n/3)^{1/4}
$$
for every planar graph $G$. Better bounds on $\fix G$ are known
for cycles~\cite{PTa}, trees~\cite{Goaoc,Bose} and, more generally, outerplanar 
graphs~\cite{SWo,RVe}. In all these cases it was shown that $\fix G=\Omega(n^{1/2})$. 

No efficient algorithm determining the parameter $\fix G$ is known.
Moreover, computing $\fix{G,\pi}$ is known to be NP-hard~\cite{Goaoc,Ver}.

\end{document}